\documentclass[amstex, 11pt]{amsart}

\input amssym.def
\input amssym.tex

\input xy
\xyoption{all}

\def\depth{\mathop{\rm depth}}
\def\grade{\mathop{\rm grade}}
\def\dim{\mathop{\rm dim}}

\def\olin{\overline}

\def\F{{\mathcal F}}

\def\m{{\mathfrak m}}

\def\Cn1{{C( x_1, {\bf x}_2, {\F}, (n_1, n_2))}}

\def\x{{\bf x}}

\newtheorem{thm}{Theorem}[section]
\newtheorem{lemma}[thm]{Lemma}
\newtheorem{cor}[thm]{Corollary}
\newtheorem{pro}[thm]{Proposition}
\newtheorem{example}[thm]{Example}
\newtheorem{remark}[thm]{Remark}
\newtheorem{defn}[thm]{Definition}
\newtheorem{notn}[thm]{Notation}
\newtheorem{blank}[thm]{}

\newcommand{\bt}{\begin{thm}}
\newcommand{\et}{\end{thm}}
\newcommand{\blem}{\begin{lemma}}
\newcommand{\elem}{\end{lemma}}
\newcommand{\bco}{\begin{cor}}
\newcommand{\eco}{\end{cor}}
\newcommand{\bp}{\begin{pro}}
\newcommand{\ep}{\end{pro}}
\newcommand{\bex}{\begin{example}}
\newcommand{\eex}{\end{example}}
\newcommand{\brm}{\begin{remark}}
\newcommand{\erm}{\end{remark}}
\newcommand{\bdefn}{\begin{defn}}
\newcommand{\edefn}{\end{defn}}
\newcommand{\bnot}{\begin{notn}}
\newcommand{\enot}{\end{notn}}
\newcommand{\bblank}{\begin{blank}}
\newcommand{\eblank}{\end{blank}}

\def\nno{\nonumber}

\newcommand{\bib}{\bibitem}
\newcommand{\f}{\frac}

\newcommand{\beqn}{\begin{eqnarray*}}
\newcommand{\eeqn}{\end{eqnarray*}}
\newcommand{\beq}{\begin{eqnarray}}
\newcommand{\eeq}{\end{eqnarray}}
\newcommand{\been}{\begin{enumerate}}
\newcommand{\eeen}{\end{enumerate}}
\newcommand{\lrar}{\longrightarrow}
\newcommand{\rar}{\rightarrow}

\newcommand{\Z}{\mathbb Z}

\newcommand{\limn}{\underset{\underset{n}{\longrightarrow}}{\lim}}
\usepackage{bbding}
\usepackage{pifont}
\usepackage{fancybox}
\usepackage{color}
\usepackage{picinpar}
\usepackage{graphicx}
\usepackage{chancery}
\usepackage{palatino}

\newcommand{\mptm}{\fontfamily{ptm}\fontsize{11}{11}\selectfont}
\newcommand{\mptmcomx}{\fontfamily{ptm}\fontsize{9}{11}\selectfont}

\newfont{\sftiny}{cmss9 scaled 1000}

\begin{document}

\title[Fiber cones]{Homology, mixed multiplicities and fiber cones}

\author{Clare D'Cruz*}
\address{Chennai Mathematical Institute, Plot H1, Sipcot IT Park, Siruseri,
Paddur PO, Siruseri Chennai 603 103, India}
\email{clare@cmi.ac.in}

\author{Anna Guerrieri}
\address
{Dipartimento di Matematica, Universit{\'a} di L'Aquila, 67100 Italy}
\email{guerran@univaq.it}




\maketitle
\mptm

\begin{abstract}
In this paper we study the Hilbert coefficients the fiber cone. We also compare the depths of the fiber cone and the associated graded ring.
\end{abstract}

\section{Introduction}

Throughout this paper we will assume that $(R,\m )$ is a local ring of
dimension $d$ with infinite residue field.  Let $I$ be an ideal in
$(R,\m)$,  then the blowing up of Spec~$R$ with center $V(I)$ is the
morphism $Proj~R[It] \lrar Spec~R$ with exceptional divisor
$Proj(R[It] \otimes_R R/I)$  and the fiber of $\m$ is $Spec~\left(R[It]
\otimes_R R/ \m \right)$.   
The ring $R(I) :=R[It] = \oplus_{n \geq 0} I^n t^n$ is Rees ring, 
$G(I) :=\oplus_{n  \geq 0} I^n/I^{n+1}$  is the associated graded ring.
 We call $F(I) := \oplus_{n   \geq 0} I^n/ \m I^n$ is the fiber cone. All these rings are called the blow-up rings.
The Rees ring and the associated graded
ring have been studied in great detail. Of recent interest is the
fiber cone. Let $I_1$ and $I_2$ be ideals in $(R,\m)$ We call  $F_{I_1}(I_2) := \oplus_{n   \geq 0} I_2^n/ I_1 I_2^n$ is the fiber cone of $I_2$ with respect to $I_1$. We are interested in the Hilbert coefficients of  $F_{I_1}(I_2)$ and its depth. 

Interest in the the Hilbert coefficients of an $\m$-primay ideal began with the work of Northcott \cite{north} and Narita \cite{narita}. Later Sally (\cite{sally1}, \cite{sally2}) realised that the depth of the $G(\m)$ played a major role in the study of the Hilbert coefficients of $\ell(R/ \m^n)$. These results were gereralised to any $\m$-primary ideals by the second author 
in her thesis \cite{anna-thesis}. She use the  Koszul complex of the Rees ring to obtain several interesting results   on the depth of the associated graded ring and its relation with the Hilbert coefficients. T.~Marley  and 
A.~Huckaba subsequently constructed a complex which they called the modified  Koszul complex via which  they  generalized the complex obtained in \cite{anna-thesis} and could get some more interesting  results \cite{tom-huc}.

The Hilbert coefficients of the fiber cone has been of interest recently. Very little is known about  the Hilbert coefficients of the fiber cone. 
In \cite{clare-verma1} and \cite{clare-verma2} we were able to obtain the Hilbert series of the fiber cone in some special cases. Later in \cite{jay-verma-jpaa}, \cite{jay-verma} and \cite{jay-tony-verma} the first Hilbert coefficient and the Hilbert series in some special cases was obtained. In all these papers it was evident  that if the associated graded ring had high depth then so did the fibercone. 
 In \cite{clare-tony}, the first author and T. Puthenpurakal, obtained some relations among the first few coefficients, under certain conditions.
In this paper we investigate the remaining coefficients.
Recently Rossi and Valla \cite{rossi-valla} obtained a few interesting results for the fibercone of a $I$-good filtration of a module.

Our approach towards the   study  of the fiber cone is  via the bigraded filtration 
${\mathcal F} = \{ I_{1}^{n_1} I_{2}^{n_2} \}_{n_1 \geq 0, n_2 \geq 0}$ where 
${\mathcal F}_1 = \{ I_{1}^{n}  \}_{n\geq 0}$  (resp. ${\mathcal F}_2 = \{ I_{2}^{n}  \}_{n\geq 0}$) is $I_1$ 
(resp. $I_2$) good filtration. 
Let 
$\x_{1k_1} = x_{11}, \ldots, x_{1k_1} \in I_{1}$ and $\x_{2k_2} = x_{21}, \ldots, x_{2k_2} \in I_{2}$. 
We   construct  a bigraded  complex  which we will denote by $C({\x}_{1k_1}, {\x}_{2k_2},{\mathcal F}, (n_1, n_2))$.
 This complex is obtained using  the mapping cone construction in a similar way as in \cite{tom-huc}. 
The beauty of this construction is that by specializing either the ideals or the exponents we are able to recover results simultaneously for the associated graded ring as well as the fiber cone. Thus we have a single complex which can give information on both the graded rings.

Using the bigraded filtration we obtain  Huneke's fundamental lemma for two ideals. 
This was obtained first by Huneke for an $\m$-primary  ideal in a two dimensional Cohen-Macaulay 
 local ring $(R, \m)$ \cite{huneke}
 and  was generalised to $\m$-primary ideals in a $d$-dimensional Cohen-Macaulay local ring in \cite{huc}. 
 In \cite{jay-verma-jpaa}, this result was generalised for the function $\ell(R/I_1 I_2^n)$, where $I_1$ and $I_2$ are $m$-primary ideals in a two dimensional Cohen-Macaulay local ring and $I_1$ is an ideal containing $I_2$. In section three we obtain  Huneke's fundamental lemma for any any two $\m$-primary ideals in a $d$-dimensional Cohen-Macaulay local ring.
 As a consequence we are able to obtain the Huneke's fundamental lemma for the fibercone $F_{I_1}(I_2)$.

  In this paper we are also able to obtain  a formula for all the Hilbert coefficients in terms of lengths of homologies of $C({\x}_{1k_1}, {\x}_{21},{\mathcal F}, (n_1, n_2))$. We give an upper bound on the first Hilbert coefficient and describe when the upper bound is attained.

 An interesting result which plays an important role in the Cohen-Macaulayness of blowup rings 
is the well known result of Goto and Shimoda which gives a necessary and sufficient condition 
for the Rees ring to be Cohen-Macaulay \cite{goto-s}. From their result it follows that if 
$(R, \m)$ is a Cohen-Macaulay local ring and the Rees ring is Cohen-Macaulay, then so is the associated graded ring. 
A natural question to ask is whether the depth of the fiber cone is 
related with the depth of the other blowup rings. In this paper we show that the vanishing of the homologies of the complex
$h_i(C({\x}_{1}, {\x}_{2k_2},{\mathcal F}, (1, n_2))$, $\depth~G(I_2)$ and $\depth~F_{I_1}(I_2)$ are inter-related.

The paper is organised as follows. In Section two describe the complex $C({\x}_{1k_1}, {\x}_{21},{\mathcal F}, (n_1, n_2))$ and develop several interesting consequences. In Section three prove the Huneke's fundamental lemma and describe the first Hilbert coefficients of the fiber cone. In section four we describe the remaining fiber coefficients.  In section five we give a relation discuss the vanishing of homologies of the complex  $C({\x}_{1}, {\x}_{2k_2},{\mathcal F}, (1, n_2))$, the depth of the fiber cone and the associated graded ring.

The results in this paper can be used to compare depths of blowup rings of contracted ideals in a regular local ring of dimension two. We can also compare the depth of the associated graded ring and ideals of minimal mixed multiplicity (\cite{clare1}, \cite{clare2}).

{\bf Acknowledgements} The first author is grateful to J. Verma for motivating the topic. The first author also wishes to thank Universit{\`a} di L'Aquila, 
for local hospitality and Progetto INDAM-GNSAGA for financial support during the preparation of this work.

\section{The Koszul Complex}

Let $(R, \m)$ be a local ring and let $I_1$ and $I_2$ be ideals of $R$.
Let ${\F}_i = \{I_{i}^{n} \}_{n \geq 0}$, $i=1,2$ be $I_{i}$-good filtrations 
and
let ${\mathcal F} = \{I_{1}^{n_1} I_{2}^{n_2} \}_{n_1,n_2 \geq 0}$. 
Let ${\bf x}_{1k}: x_{11}, \ldots, x_{1k_{1}} \in I_1$ and let 
${\bf x}_{2k}: x_{21}, \ldots, x_{2k_{2}} \in I_2$. We construct the 
bigraded complex:
$C({\bf x}_{1k_1}, {\bf x}_{2k_2}, {\mathcal F}, (n_1, n_2))$ as follows:

Using the concept of mapping cylinder we can construct the following complexes (\cite {tom-huc}):
\mptmcomx
\beq
\label{newdisplay}
 C({\bf x}_{1k_1}, {\mathcal F_1}, (n_1, n_2)): 
      0
\rar \f{R}{I_1^{n_1-k_1} I_2^{n_2} }
\rar \cdots 
\rar \left( \f{R}{I_1^{n_1-i} I_2^{n_2}} \right)^{k_1 \choose i}
\rar \cdots 
\rar \left( \f{R}{I_1^{n_1-1} I_2^{n_2}} \right)^{k_1 \choose 1}
\rar \f{R}{I_1^{n_1} I_2^{n_2}} \rar 0
\eeq

\beq
\label{newdisplay-2}
 C({\bf x}_{2k_2}, {\mathcal F_2}, (n_1, n_2)):     
 0
\rar \f{R}{I_1^{n_1} I_2^{n_2- k_2} }
\rar \cdots 
\rar \left( \f{R}{I_1^{n_1} I_2^{n_2-i}} \right)^{k_2 \choose i}
\rar \cdots 
\rar \left( \f{R}{I_1^{n_1} I_2^{n_2-1}} \right)^{k_2 \choose 1}
\rar \f{R}{I_1^{n_1} I_2^{n_2}} \rar 0
\eeq

\mptm
We also have the mapping of complexes:
${.x_{21}}: C({\bf x}_{1k_1}, {\mathcal F}, (n_1, n_2-1)) \lrar C({\bf x}_{1k_1}, {\mathcal F}, (n_1, n_2))$

\beqn
\xymatrix@C=10pt
{
0 \ar[r]
&  \f{R}{I_1^{n_1-k_1} I_2^{n_2-1}}  \ar[d]_{.x_{21}}   \ar[r]
&  \left( \f{R}{I_1^{n_1-k_1+1} I_2^{n_2-1}} \right)^{k_1 \choose k_1 -1}  \ar[d]_{.x_{21}}   \ar[r]
& \cdots \ar[r]
&      \left(  \f{R}{I_1^{n_1-2} I_2^{n_2-1}} \right)^{k_1 \choose 2}  \ar[d]_{.x_{21}} \ar[r]
&       \left( \f{R}{I_1^{n_1-1} I_2^{n_2-1}} \right)^{k_1 \choose 1}  \ar[d]_{.x_{21}} \ar[r]
&  \f{R}{I_1^{n_1} I_2^{n_2-1}}  \ar[d]_{.x_{21}}\ar[r] 
& 0\\
0 \ar[r]
& \f{R}{I_1^{n_1-k_1} I_2^{n_2}}  \ar[r]
& \left( \f{R}{I_1^{n_1-k_1+1} I_2^{n_2}} \right)^{k_1 \choose k_1 -1} \ar[r]
& \cdots \ar[r]
&       \left( \f{R}{I_1^{n_1-2} I_2^{n_2}} \right)^{k_1 \choose 2}   \ar[r]
&       \left(  \f{R}{I_1^{n_1-1} I_2^{n_2}} \right)^{k_1 \choose 1}   \ar[r]
&  \f{R}{I_1^{n_1} I_2^{n_2}}  \ar[r] 
& 0
}
\eeqn

The  mapping cylinder of these two complexes is the following complex which we denote by  $C({\bf x}_{1k_1}, {x}_{21},{\mathcal F}, (n_1, n_2))$ and is of the form:
\beqn
\xymatrix@C=10pt
{    0  
\rar  \f{R}{I_1^{n_1-k_1} I_2^{n_2-1}} 
\rar  \cdots 
\rar    
       \left( 
       \f{R}
        {I_1^{n_1-i+1} I_2^{n_2-1}} 
       \right)^{k_1 \choose i-1}
       \oplus 
       \left(   
       \f{R}
         {I_1^{n_1-i} I_2^{n_2}} 
       \right)^{k_1 \choose i}
\rar   \cdots   
\rar   \f{R}
         {I_1^{n_1} I_2^{n_2-1}}  
       \oplus 
       \left(   
       \f{R}
         {I_1^{n_1-1} I_2^{n_2}} 
       \right)^{k_1 \choose 1}  
\rar  \f{R}
        {I_1^{n_1} I_2^{n_2}}  
\rar     0. }
\eeqn

Inductively, for $k_1, k_2 \geq 0$ we get the complex which we denote by $ C({\x}_{1k_1}, {\x}_{2k_2},{\mathcal F}, (n_1, n_2))$:
\beq
\label{main-complex}
\xymatrix@C=10pt{   0 
\rar  \f{R}{I_1^{n_1-k_1} I_2^{n_2-k_2}}     
\rar  \cdots  
\rar  \bigoplus_{j=0}^{i} 
      \left( 
      \f{R}
        {I_1^{n_1-i+j} I_2^{n_2-j}} 
       \right)^{{k_1 \choose i-j}
               {k_2 \choose j}}   
\rar   \cdots  
\rar   \f{R}{I_1^{n_1} I_2^{n_2}}  
\rar    0.}
\eeq

\brm
Note that if $k_1 = 0$, then,
\beqn
 C({\bf x}_{1k_1}, {\bf x}_{2k_2},{\mathcal F}, (n_1, n_2))
=  C(                {\bf x}_{2k_2},{\mathcal F}, (n_1, n_2)) .
\eeqn
\erm

We have the short exact sequence of complexes:
\mptmcomx
\beq
\label{ses-comx}
     0 
\rar {{C({\bf x}_{1k_1},  {\bf x}_{2,k_2 -1},  {\mathcal F}, (n_1, n_2))}}
\rar {{C({\bf  x}_{1k_1}, {\bf  x}_{2, k_2},  {\mathcal F}, (n_1, n_2))}} 
\rar {{C({\bf x}_{1k_1},  {\bf x}_{2,k_2 -1}, {\mathcal F}, (n_1, n_2-1))}}[-1]
\rar 0
\eeq
\mptm
and the   sequence of homologies:
\mptmcomx
\beq
\label{main-homology}
          \cdots
\rar &  H_i( {{C({\bf  x}_{1  k_1}, 
                  {\bf  x}_{2  k_2},  
                  {\mathcal F}, (n_1, n_2))}} )
\rar &    H_{i-1}({{C( {\bf x}_{1 k_1}, 
                       {\bf x}_{2 k_2 -1},  
                       {\mathcal F}, (n_1, n_2-1))}})
\\ \nno
\rar      H_{i-1}({{C( {\bf x}_{1k_1}, 
                       {\bf x}_{2k_2 -1},  
                       {\mathcal F}, (n_1, n_2))}})
\rar&    H_{i-1}( {{C({\bf  x}_{1  k_1}, 
                       {\bf  x}_{2  k_2},  
                       {\mathcal F}, (n_1, n_2))}} )
\rar& \cdots
\eeq
\mptm

\bt
\label{homology}
Let $(R,\m)$ be a Cohen-Macaulay local ring and let $I_1$ and $I_2$ be $\m$-primary ideals of $R$.
Let  $k_1, k_2 \geq 0$ with $k_1 + k_2 \geq 1$.
 Let $\x_{1k_1} \in I_1$ 
 and $\x_{2k_2} \in I_2$ be a superficial sequence
for $I_1$ and $I_2$ which is a regular sequence. Then
\been
\item
\label{homology-one}
For all $n_1, n_2 \in \Z$, 
${\displaystyle
        H_0( {{C({\bf  x}_{1  k_1}, 
                  {\bf  x}_{2  k_2},  
                  {\mathcal F}, (n_1, n_2))}} )
                  \cong \f{R}{I_1^{n_1} I_2^{n_2} + (\x_{1, k_1}, \x_{2k_2})}.
}$\\

\item
\label{homology-two}
For all $n_1 \geq 1$ and $n_2 \geq 1$, 
${\displaystyle
    H_1( {{C({\bf  x}_{1  k_1}, 
                  {\bf  x}_{2  k_2},  
                  {\mathcal F}, (n_1, n_2))}} )
\cong \f{(\x_{1 k_1}, \x_{2 k_2}) \cap  I_1^{n_1} I_2^{n_2}}
  { \x_{1 k_1} I_1^{n_1-1} I_2^{n_2} + \x_{2 k_2} I_1^{n_1} I_2^{n_2-1}}}.
$\\

\item
\label{homology-three}
For all $n_1, n_2 \in \Z$, 
${\displaystyle
     H_{k_1 + k_2}( {{C({\bf  x}_{1  k_1}, 
                  {\bf  x}_{2  k_2},  
                  {\mathcal F}, (n_1, n_2))}} )
\cong \f{  (I_1^{n_1- k_1+1} :(\x_{1 k_1}) )
          \cap (I_2^{n_2-k_2+1} : (\x_{2 k_2})) }
        {I_1^{n_1-k_1} I_2^{n_2-k_2}}}.$

\eeen
\et
\proof
(1) and (3) are easy to verify. We prove (2).
Consider
\beqn
      \psi_{1} 
& =&  \left( \f{R}{I_1^{n_1-1} I_2^{n_2}}   \right)^{k_1 \choose 1}       
      \oplus 
      \left( \f{R}{I_1^{n_1}   I_2^{n_2-1}} \right)^{k_2 \choose 1}  
 {\rar}  
         \f{R}{I_1^{n_1} I_2^{n_2}} 
\rar      0 \\
\psi_2 &=& 
       \left( 
       \f{R}
         {I_1^{n_1-2} I_2^{n_2}} \right)^{k_1 \choose 2}
\oplus  \left( 
       \f{R}
         {I_1^{n_1-1} I_2^{n_2-1}} 
       \right)^{{k_1 \choose 1}{k_2 \choose 1}} 
\oplus \f{R}
         {I_1^{n_1} I_2^{n_2-2}}^{k_2 \choose 2} 
\lrar  \left( 
       \f{R}
         {I_1^{n_1-1} I_2^{n_2}}   \right)^{k_1 \choose 1}       
      \oplus 
      \left( 
      \f{R}
        {I_1^{n_1}   I_2^{n_2-1}} \right)^{k_2 \choose 1}   
\eeqn
We claim that 
\beqn
      \f{\mbox{Ker~} \psi_{1}}{\mbox{Im} \psi_{2}} 
\cong \f{(\x_{1 k_1}, \x_{2 k_2}) \cap  I_1^{n_1} I_2^{n_2}}
        { \x_{1 k_1} I_1^{n_1 -1} I_2^{n_2} 
        + \x_{2 k_2} I_1^{n_1} I_2^{n_2}}.
\eeqn
Put $x_i = x_{1i}, 1 \leq i \leq k_1$ and $x_{j+i} = x_{2j}, 1 \leq j \leq k_2$.
Consider that map
\beqn
   \phi: (Ker\, \psi_1) 
\, \lrar 
\f{(x_1, \ldots, x_{k_1 + k_2}) \cap  I_1^{n_1} I_2^{n_2}}
  { (x_{1} \ldots  x_{k_1})I_1^{n_1 -1} I_2^{n_2} 
  + (x_{k_1 + 1}, \ldots, x_{k_1 + k_2}) I_1^{n_1} I_2^{n_2}}   
      \eeqn
defined by
$
      \phi(a_1', \ldots,  a_{k_1 + k_2}') 
=   (x_{1} a_1 + \cdots + x_{k_1+k_2} a_{k_1+k_2})'$ where 
primes  denote the respective residue classes. 

{\bf Claim:} There  exist elements 
$b_1, \ldots, b_{k_1}  \in I_1^{n_1-1}I_2^{n_1}$, 
$b_{k_1 + 1}, \ldots, b_{k_1 + k_2} \in I_1^{n_1} I_2^{n_2-1}$ 
and elements $r_{ij}, r_{ji} \in R$ 
such
that 
\beq
\label{eq1}
a_i = b_i - \sum_{j=1}^{i-1} x_j r_{ji}
          + \sum_{j=i+1}^{k_1 + k_2} x_j r_{ij}
\eeq
and equations
\beq
\label{eq2}
\sum_{k= i+1}^{k_1 + k_2} x_k
\left[ 
a_{k} +  \sum_{j=1}^{i-1} x_j r_{jk}
          \right]
= \sum_{k=i+1}^{k_1 + k_2} x_k b_k 
\eeq
where the terms for which the subscripts are negative are zero. 

We prove by induction on $i$.
Let i = 1. 
 Since  
 $(a_1', \ldots, a_{k_1 + k_2}') \in Ker \, \psi_1$,
there exist elements 
$b_{1}, \ldots, b_{k_1} \in I_1^{n_1-1} I_2^{n_2}$ 
and elements 
$b_{k_1 + 1},  \ldots, b_{k_1 + k_{2}}$
such that 
\beq
\label{eq3} 
     x_{1} a_{1} + \cdots + x_{k_1+ k_2} a_{k_1+ k_2}
&=&  x_{1} b_{1} + \cdots + x_{ k_1+ k_2} b_{k_1+k_2}.
\eeq
 Then
\beqn
        x_{1}(a_{1} - b_{1}) \in (x_{2}, \ldots,  x_{k_1+k_2}).
\eeqn
Since 
$x_1, \ldots, x_{k_1+ k_2}$  
is a regular sequence in $R$ there exist elements $r_{1j} \in R$ 
such that 
\beq
\label{eq4}
a_{1} = b_{1} + \sum_{j=2}^{k_1+k_2} x_{j} r_{1j} .
\eeq
Substituting (\ref{eq4}) in (\ref{eq3}) we get
\beq
\label{eq5} 
      x_{1} \left[b_{1} + \sum_{j=2}^{k_1+k_2} x_{j} r_{1j} \right] 
    + x_{2} a_{2} + \ldots + x_{k_1 + k_2} a_{k_1 + k_2}
&=&   x_{1} b_{1} + \ldots + x_{k_1 + k_2} b_{k_1 + k_2}.
\eeq
which gives
\beqn
  \sum_{j=2}^{k_1 + k_2} x_{j}
  \left[ x_{1}  r_{1j} 
+  x_{j} a_{j} \right]
= \sum_{j=2}^{k_1 + k_2} x_{j} b_{j}.
\eeqn
This proves the claim for $i=1$. 

Now suppose $i > 1$ and that the claim is true for i. Then by induction hypothesis we have
\beq
\label{eq6}
\sum_{k= i+1}^{k_1 + k_2} x_k
\left[ 
a_{k} +  \sum_{j=1}^{i} x_j r_{jk}
          \right]
= \sum_{k=i+1}^{k_1 + k_2} x_k b_k .
\eeq
Hence
\beq
\label{eq7}
       x_{i+1}
   \left[ a_{i+1} - b_{i+1} 
+  \sum_{j=1}^{i} x_j r_{j,i+1}\right]
+  \sum_{k= i+2}^{k_1 + k_2} x_k
   \left[ 
    a_{k} 
+  \sum_{j=1}^{k-1} x_j r_{jk}
          \right]
= \sum_{k=i+2}^{k_1 + k_2} x_k b_k .
\eeq
Therefore
\beq
\label{eq8}
         x_{i+1}
   \left[ a_{i+1} - b_{i+1} 
+  \sum_{j=1}^{i} x_j r_{j,i+1}\right]
\in (x_{i+2}, \ldots, x_{k_1 + k_2}).
\eeq
Since $x_{i+1}, \ldots, x_{k_1 + k_2}$ is a regular sequence,
there exists elements $r_{i+1,j} \in R$, $i+2 \leq j \leq k_1 + k_2$ such that 
 \beq
\label{eq9}
 a_{i+1} - b_{i+1} +  \sum_{j=1}^{i} x_j r_{j,i+1}
= \sum_{j=i+2}^{k_1 + k_2} x_j r_{i+1,j}.
\eeq
Substituting (\ref{eq9}) in (\ref{eq7}) we get 
\beq
\label{eq10}
        x_{i+1}
     \left[ 
     \sum_{j=i+2}^{k_1 + k_2} x_j r_{i+1,j}\right]
+    \sum_{k=i+2}^{k_1 + k_2} x_k
     \left[ a_{k} 
+  \sum_{j=1}^{i} x_j r_{jk}          \right]
= \sum_{k=i+2}^{k_1 + k_2} x_k b_k .
\eeq
The left side of the above equation is:
\beq
\label{eq10}
\sum_{k= i+2}^{k_1 + k_2}
x_k \left[ a_k + \sum_{j=1}^{i} x_j r_{jk} \right].
\eeq
This completes the proof of the theorem.
\qed

We now describe the vanishing of some of these homologies.

\blem
\label{vanishing}
Let $(R, \m)$ be a Cohen-Macaulay local ring. Let $I_1$ and $I_2$ be $\m$-primary ideals of $R$. Let $k_1$ and $k_2$ be non-negative integers with $k_1 + k_2 \geq 1$. Let $x_{11}, \ldots, x_{1k_1} \in I_1$  and $x_{21}, \ldots, x_{2k_2} \in I_2$.
 such that $\x_{1k_1}, \x_{2k_2}$ is a regular sequence which is superficial for $I_1$ and $I_2$. 
\been
\item
\label{vanishing-one}
Let $k_1 = 0$  Then 
 for all $i \geq 1$ and for all $n_2 \gg 0$, 
 $H_i(C({\x}_{2k_2}, {\mathcal F}, (n_1, n_2)))=0$. 

\item
\label{vanishing-two}
Let $k_1 = 1$ and $k_2 \geq 1$. Then
for all $i \geq 2$   and for all $n_2 \gg 0$, 
 $
H_i(C(x_1, \x_{2k_2}, {\mathcal F}, (n_1, n_2)))
=0$.

 \item
 \label{vanishing-three}
Let $k_1,k_2 \geq 2$.
For all $i \geq 1$,  
$H_i( C(\x_{1~k_1}, \x_{2~k_2}, {\mathcal F}, (n_1, n_2))) = 0$  
 if  $n_1\gg 0$ and  $n_2 \gg 0$.
\eeen
\elem
\proof Let $k_1 = 0$ and $k_2 = 1$. Then by Theorem~\ref{homology}(\ref{homology-three})
\beqn
            H_1( {{C(x_1,  {\mathcal F}, (n_1, n_2))}})
&\cong& \f{  I_1^{n_1} I_2^{n_2} : (x_1) }{I_1^{n_1} I_2^{n_2-1}}.
\eeqn
Since $x_1$ is regular and superficial for ${\mathcal F}$, $n_2 \gg 0$, 
 $I_1^{n_1} I_2^{n_2} : (x_{1}) = I_1^{n_1} I_2^{n_2-1}$. 
Hence
$H_1( {{C(x_1,  {\mathcal F}, (n_1, n_2))}})  = 0$ for all $n_2 \gg 0$.

Let  $k_2 > 1$ and assume that (1) is true for $k_1-1$. We have the exact sequence of complexes
\beqn
0 \rar {{C( {\bf x}_{2~k_2 -1},  {\mathcal F}, (n_1, n_2))}}
 \rar
 {{C({\bf  x}_{2 ~ k_2},  {\mathcal F}, (n_1, n_2))}} 
   \rar
   {{C( {\bf x}_{2 ~k_2 -1},  {\mathcal F}, (n_1, n_2-1))}}[-1]
   \rar 0 
\eeqn
gives rise to the   sequence:
\beq
\label{complex-1}
\begin{array}{rlllll}
 \cdots 
&\rar&    H_i( {{C({\bf  x}_{2k_2},  {\mathcal F}, (n_1, n_2))}} )
&\rar&    H_{i-1}({{C( {\bf x}_{2~k_2 -1},  {\mathcal F}, (n_1, n_2-1))}})
\\ 
\rar           H_{i-1}({{C( {\bf x}_{2~k_2 -1},  {\mathcal F}, (n_1, n_2))}})
&\rar&    H_{i-1}( {{C({\bf  x}_{2~k_2-1},  {\mathcal F}, (n_1, n_2))}} )
&\rar& \cdots 
\end{array}
\eeq
Put $i=1$ in the complex (\ref{complex-1}). For $n_2 \gg 0$ we have the exact sequence
\mptmcomx
\beq
\label{ses-new}
       0 
\rar H_1( {{C({\bf  x}_{2 k_2},  {\mathcal F}, (n_1, n_2))}} )
\rar H_{0}({{C( {\bf x}_{2k_2 -1},  {\mathcal F}, (n_1-1, n_2))}}) 
\rar  
     H_{0}({{C( {\bf x}_{2k_2 -1},  {\mathcal F}, (n_1, n_2))}}
\rar 0
\eeq
\mptm
since $H_1( {{C({\bf  x}_{2 k_2-1},  {\mathcal F}, (n_1, n_2))}} )=0$ for 
all $n_2 \geq 0$.

By Theorem~\ref{homology} and the sequence (\ref{ses-new}) we have
\beqn
0  \lrar H_1 (  {{C( {\bf x}_{2 ~k_1},  {\mathcal F}, (n_1, n_2))}})
 \lrar  \f{R}{I_1^{n_1} I_2^{n_2-1} + ({\bf  x}_{2~k_2 -1})}
   \lrar   \f{R}{I_1^{n_1} I_2^{n_2} + ({\bf  x}_{2~k_2 -1})}
   \lrar 0 .
   \eeqn
   
   Let $\olin{\mathcal F}$ denote the filtration 
   ${\mathcal F}/ ({\bf x}_{2~k_2 -1})$. Let 
   $\olin{x_{2~k_2}}$ denote the image of $x_{2~k_2}$ in 
   $R/( {\bf  x}_{2~k_2-1})$. Then 
   by the case $k_2=1$, 
$H_1( C( \olin{x_{2~k_2}}, \olin{\mathcal F}, (n_1, n_2))) = 0$ for all $n_2 \gg 0$. 
If $i > 1$, then 
by induction on $k_2$, for all $i >1$, $H_i( {{C({\bf  x}_{2 ~ k_2-1},  {\mathcal F}, (n_1, n_2))}} ) = 0$ for all $n_2 \gg 0 $. Hence from (\ref{complex-1}), $H_i( {{C({\bf  x}_{2 ~ k_2-1},  {\mathcal F}, (n_1, n_2))}} ) = 0$ for all $n_2 \gg 0 $.
This proves (1). 

Now let $k_1 = 1$. 
Now if $i = 1$ and if $n_1 \gg 0$, then once again by the argument similar to (1),
we can show that  for all $n_1, n_2 \gg 0$, 
$H_i ({{C(\x_{1k_1} x_{21},  {\mathcal F}, (n_1, n_2))}} )= 0$.
Let $i >1$. 
Consider
\beqn
       0 
\lrar {{C(\x_{2k_2},  {\mathcal F}, (n_1, n_2))}}
\lrar {{C(x_1,\x_{2k_2} ,  {\mathcal F}, (n_1, n_2))}} 
\lrar {{C(\x_{2k_2},  {\mathcal F}, (n_1-1, n_2))}}[-1]
\lrar 0. 
\eeqn

This gives rise to the exact sequence
\beqn
      \cdots 
\rar H_i ({{C(\x_{2k_2},  {\mathcal F}, (n_1, n_2))}})
\rar H_i ({{C(x_1, \x_{2k_2}  {\mathcal F}, (n_1, n_2))}} )
\rar H_{i-1}( {{C(\x_{2k_2},  {\mathcal F}, (n_1-1, n_2))}})
\rar \cdots
\eeqn
By the previous case
$H_i ({{C(\x_{2k_2},  {\mathcal F}, (n_1, n_2))}})
=  H_{i-1}( {{C(\x_{2k_2},  {\mathcal F}, (n_1-1, n_2))}}) =0$
for all $i \geq 2$ and for all $n_2 \gg 0$. This proves (2).

   (3) follows by induction on $k_1$. 
   \qed

         \bnot
   \beqn
     h_i (\x_{1 k_1}, \x_{2 k_2})(n_1,n_2) 
&=& \ell
    \left( H_i (\x_{1 k_1}, \x_{2 k_2}, {\mathcal F},(n_1,n_2)) 
    \right)\\
     h_i(\x_{1k_1,}, \x_{2 k_2})(n_1, *) 
&=& \sum_{n_2 \geq 0} 
     h_i (\x_{1 k_1}, \x_{2 k_2})(n_1,n_2) \\
     h_i(\x_{1k_1}, \x_{2 k_2})(*, n_2) 
&=& \sum_{n_1 \geq 0} 
     h_i (\x_{1 k_1}, \x_{2 k_2})(n_1,n_2) \\
   h_i(\x_{1k_1}, \x_{2 k_2}) 
&=&    \sum_{n_1 \geq 0} h_i(\x_{1k_1}. \x_{2 k_2})(n_1) \\
&=& \sum_{n_2 \geq 0}h_i(\x_{1k_1}, \x_{2 k_2})(n_2) \\
 h_i (\x_{1 k_1})(n_1,n_2) 
&=& \ell
    \left( H_i (\x_{1 k_1},  {\mathcal F},(n_1,n_2)) 
    \right)\\
     h_i ( \x_{2 k_2})(n_1,n_2) 
&=& \ell
    \left( H_i ( \x_{2 k_2}, {\mathcal F},(n_1,n_2)) 
    \right)\\
            \eeqn 
               \enot
               
               The complex $\Cn1$ satisfies a certain rigidity. The rigidity condition was proved for  the complex
               $C(\x_{2k_2},  {\mathcal F}, (0, n_2))$ was proved in \cite[Lemma~3.4]{tom-huc}. The same proof can be adopted to show the rigidity of  the complex $C(\x_{2k_2},  {\mathcal F}, (1, n_2))$ and 
               $C(x_1, \x_{2k_2},  {\mathcal F}, (1, n_2))$.
               
                 \blem
\label{rig-0}
Let $i \geq 1$ and $k_2 \geq 1$.    Suppose 
   $H_i ({{C( \x_{2k_2},  {\mathcal F}, (1, n_2))}} )=0$ for all $n_2 \geq 0$. Then $H_j ({{C( \x_{2k_2},  {\mathcal F}, (1, n_2))}} )=0$ for 
   all $j > i$ and $n_2 \geq 0$.
      \elem
\proof Let $k_2 = 1$, there is nothing to prove. Hence we can assume that $k_2 \geq 2$. The proof here is similar to that of \cite[Lemma~3.4]{tom-huc} applied to the complex
\beq
   \label{complex-rigidity-2-1}
\begin{array}{rrlllll}
     &    0  
&\rar&    H_{2}( {{C(\x_{2k_2},   {\mathcal F}, (1, n_2))}} )
&\rar&    H_{1}({{C(x_{2k_2-1},        {\mathcal F}, (1, n_2-1))}})
\\ 
\rar&     H_{1}({{C(x_{2k_2-1},         {\mathcal F}, (1, n_2))}})
&\rar&    H_{1}( {{C(\x_{2k_2},   {\mathcal F}, (1, n_2))}} )
&\rar&    \cdots. \\
\end{array}
\eeq
  \qed

\blem
\label{rig-1}
Let $i, k_2 \geq 1$.    Suppose 
   $H_i ({{C(x_1, \x_{2k_2},  {\mathcal F}, (1, n_2))}} )=0$ for all $n_2 \geq 0$. Then $H_j ({{C(x_{1}, \x_{2k_2},  {\mathcal F}, (1, n_2))}} )=0$ for 
   all $j > i$ and $n_2 \geq 1$.
      \elem
  \proof 
  Suppose 
\beqn
H_i ({{C(x_1, \x_{2k_2},  {\mathcal F}, (1, n_2))}} )=0
\hspace{.5in} \mbox{ for all } n_2 \geq 0,
\eeqn 
for some $i \geq 1$. We have the exact sequence
\beq
\label{complex-rigidity-3}
\begin{array}{rrlllll}
&      &                   \cdots  
&\rar & H_{i+1}( {{C(x_1, \x_{2k_2},   {\mathcal F}, (1, n_2))}} ) 
&\rar&   H_{i}  ( {{C(\x_{2k_2},{\mathcal F}, (1, n_2-1))}})\\ 
&\rar&  H_{i}  ( {{C(\x_{2k_2},         {\mathcal F}, (1, n_2))}})  
&\rar&    H_{i}( {{C(x_1, \x_{2k_2},   {\mathcal F}, (1, n_2))}} )=0   
&    & 
\end{array}
\eeq
How use Lemma~\ref{rig-0} and the proof is similar to \cite[Lemma~3.4]{tom-huc}.

\qed

  
  The following lemma was proved for the complex $C( \x_{2k_2},   {\mathcal F}, (0, n_2))$. Replacing by the above complex $C( x_1, \x_{2k_2},   {\mathcal F}, (1, n_2))$ (resp. $C(x_1, \x_{2k_2},   {\mathcal F}, (1, n_2))$)  we get Lemma~\ref{rigidity-new} (resp, \ref{rigidity}).
%

 \blem
   \label{rigidity-new}
    Let $i \geq 1$ and $k_2 \geq 2$. Then  
    \beqn
 \sum_{j \geq i} (-1)^{j-i} h_j(\x_{2 k_2})(1, *) \geq 0.
\eeqn
Equality holds if and only if
for all  $j \geq i$, $n_2 \geq 0$ and $0 \leq a \leq k$ we have
\beqn
 H_{j}( {{C( \x_{2k_2-a}, {\mathcal F}, (1, n_2))}})=0. 
\eeqn 
\elem
\proof The proof is similar to the proof of \cite[Theorem~3.7]{tom-huc}.
\qed

   \bt
   \label{rigidity}
    Let $i \geq 1$. Then  
    \beqn
 \sum_{j \geq i} (-1)^{j-i} h_j(x_{1}, \x_{2 k_2})(1, *) \geq 0.
\eeqn
Equality holds if and only if
$j \geq i$, $n_2 \geq 0$ and $0 \leq a \leq k$ we have
\beqn
 H_{j}( {{C( {x}_{1}, \x_{2k_2-a}, {\mathcal F}, (1, n_2))}})=0. 
\eeqn 
\et
\proof The proof is similar to the proof of \cite[Theorem~3.7]{tom-huc}.
\qed


\section{Huneke's fundamental lemma and multiplicity of the fiber cone}    

Let $(R, \m)$ be a local ring and let $I$ be an $\m$-primary ideal of $R$.  It is well known that for $n \gg 0$, the function $H(n):= \ell(R/I^n)$ is a polynomial which we will denote by $P(n)$. For a two dimensional Cohen-Macaulay local ring and an $\m$-primary ideal $I$, Huneke gave a relation between $\Delta^2[P(n) - H(n)]$ and the multiplicity of the ideal $I$ which is known as 
 Huneke's fundamental (\cite[Lemma~2.4]{huneke}). This was generalised for an $\m$-primary ideal in a d-dimensional Cohen-Macaulay local ring in \cite{huc} and for an Hilbert filtration in \cite{tom-huc}. In \cite{tom-huc} they showed that the function  $\Delta^d [P(n) = H(n)]$ and the lengths of the length of homology modules namely $h_i ( \x_{2 k_2})(0,n_2)$. In  \cite[Corollary~3.2]{jay-verma-jpaa}, Huneke's fundamental lemma was generalised for the fucntion $\ell(R/KI^n)$, where $I$ is an $\m$-primary ideal in a two dimensional Cohen-Macaulay local ring and  $K$ is an ideal containing $I$.

We state the most general form of the Huneke's fundamental lemma for a $d$-dimensional Cohen-Macaulay local ring. This generalizes the earlier results. 
 
 \bnot
Let $f: \Z^2 \lrar \Z^2$ be a function. Define
\beqn
\Delta^{k_1,k_2}f(n) :=
 \sum_{i=0}^{k_1 + k_2} (-1)^i
    \sum_{j=0}^{i} {k_1 \choose i-j} {k_2 \choose j}
    f(n_1 - (i-j),n_2-j).
    \eeqn

 \enot
\bt
[Huneke's fundamental lemma for two ideals]
\label{extension-huneke}
Let $(R, \m)$ be a Cohen-Macaulay local ring and let  $I_1$ and $I_2$ be $\m$-primary ideals and $k_1$ and $k_2$ 
non-negative integers with $k_1 + k_2 \geq 1$. Let $\x_{1k_1} \in I_1$ and $\x_{2k_2}\in I_2$. Assume that $\x_{1 k_1}, \x_{2k_2}$ is a Rees-superficial sequence for $I_1$ and $I_2$. Then
\beqn
\label{extension-huneke-general}
&& \Delta^{k_1,k_2}    H_{\mathcal F}(n_1,n_2)
    \\ \nno
&=&   \ell 
    \left( 
    \f{R}{I_1^{n_2} I_2^{n_2} + (\x_{1k_1},\x_{2 k_2})} 
    \right)
-   \ell 
    \left( 
    \f{(\x_{1 k_2}, \x_{2 k_2}) \cap I_1^{n_1} I_2^{n_2}}
       {\x_1 I_1^{n_1 -1} I_2^{n_2} + \x_2 I_1^{n_1} I_2^{n_2-1}}
    \right) 
                + \sum_{i \geq 2} (-1)^i 
                 h_i(\x_{1 k_1}, \x_{2 k_2})(n_1, n_2).
                 \eeqn
                               \et
                               \proof From Theorem~\ref{homology} and the complex (\ref{main-complex}) we have
                               \beqn
                                &&  \sum_{i=0}^{k_1 + k_2} (-1)^i
    \sum_{j=0}^{i} {k_1 \choose i-j} {k_2 \choose j}
    H_{\mathcal F}(n_1 - (i-j),n_2-j)
       \\
&=& \sum_{i \geq 0} (-1)^i 
                 h_i(\x_{1 k_1}, \x_{2 k_2})(n_1, n_2)\\
                 &=& \ell \left(
                 \f{R}{I_1^{n_1} I_2^{n_2} + (\x_{1, k_1}, \x_{2k_2})} 
                 \right)
                 - \ell \left( 
                 \f{(\x_{1 k_1}, \x_{2 k_2}) \cap  I_1^{n_1} I_2^{n_2}}
  { \x_{1 k_1} I_1^{n_1-1} I_2^{n_2} + \x_{2 k_2} I_1^{n_1} I_2^{n_2-1}}
  \right)
  + \sum_{i \geq 2} (-1)^i 
                 h_i(\x_{1 k_1}, \x_{2 k_2})(n_1, n_2).
                                   \eeqn             \qed
   
  Let $I_1$ and $I_2$ be $\m$-primary ideals and let ${\mathcal F} = I_1^{n_1}I_2^{n_2}$. It was proved in \cite{b} that  for
large values of $r$ and $s,$ the
function $  H_{\mathcal F}(n_1, n_2):=\ell(R/I_1^{n_1}I_2^{n_2})$ is given by a polynomial $P_{\mathcal F}(n_1, n_2)$ of total
degree $d$ in $r$ and $s$ and can be written in the form:
\beqn
 P_{\mathcal F}(n_1, n_2)=\sum_{i_1+i_2 \leq d}e_{i_1i_2}{n_1+i_1 \choose i_1}{n_1+i_2 \choose i_2},\eeqn
where $e_{i_1i_2}$ are certain integers. When $i_1+i_2=d,$ we put
$e_{i_1i_2}=e_{i_2}(I_1|I_2)$
for $i_2=0,1, \ldots,d.$ In this case these are called the mixed
multiplicities of $I_1$ and $I_2$. Rees introduced superficial sequence, what is well known as the Rees superficial sequence in order to explain mixed multiplicities. An element $x_1 \in I_1$ is superficial for $I_1$ and $I_2$  if  for all $n_1 \gg 0$ and for all $n_2 \geq 0$, we have
$x_1 \cap I_1^{r_1} I_2^{r_2} = x_1 I_1^{n_1-1} I_2^{n_1}$. A sequence of elements $x_{1}, x_{2}, \ldots, x_{k_1},x_{k_1+1}, x_{k_1+2}, \ldots, x_{k_1+k_2}$ where
$x_{1}, x_{12}, \ldots, x_{k_1}\in I_1$ and $x_{k_1+1}, x_{k_1+2}, \ldots, x_{k_1+ k_2} \in I_2$ is a superficial sequence if $\olin{x_{i}}$ is superficial for 
$\olin{I_1}$ and $\olin{I_2}$ where $\olin{\hphantom{X}}$ denotes the image modulo $x_1, \ldots, x_{i-1}$.

It  is well known that $e_i(I_1|I_2)$ is the multiplicity of a superficial sequence $\x_{1,d-i},\x_{2,i}$ of $I_1$ and $I_2$. In particular, $e_0(I_1|I_2)= e(I_1)$ and
$e_d(I_1|I_2)= e(I)$ where $e(.)$ denotes the Hilbert-Samuel
multiplicity \cite{r1}.  The other mixed multiplicities too can be shown
to be
Hilbert-Samuel multiplicities of certain systems of parameters (\cite{t}
and \cite{r2}). 

 
                               
     \bco
[Huneke's fundamental lemma and the difference function]
\label{cor-extension-huneke-d}
Let $(R, \m)$ be a Cohen-Macaulay local ring and let  $I_1$ and $I_2$ be $\m$-primary ideals and $k$ and non-negative integer.  Let $\x_{1k} \in I_1$ and $\x_{2.d-k}\in I_2$. Assume that $\x_{1, k}, \x_{2,d-k}$ is a Rees-superficial sequence for $I_1$ and $I_2$. Then
\beq
\label{extension-huneke-general}
 && \Delta^{k,d-k} 
  \left[ P_{\mathcal F}(n_1,n_2)
  - H_{\mathcal F}(n_1,n_2) \right] \\ \nno
    &=&    \ell 
    \left( 
    \f{ I_1^{n_1} I_2^{n_2}}
       {\x_{1k} I_1^{n_1 -1} I_2^{n_2} + \x_{2,d-k} I_1^{n_1} I_2^{n_2-1}}
    \right) 
                -\sum_{i \geq 2} (-1)^i 
                 h_i(\x_{1 k}, \x_{2,d- k})(n_1, n_2).
                 \eeq
\eco
\proof Put $k_1 = k$ and $k_2 = d-k$ in Theorem~\ref{extension-huneke}.
the right hand side of equation (\ref{extension-huneke-d}) is 
\beq
 \label{extension-fiber-d-2}
 && \ell \left(
         \f{R}{I_1^{n_1} I_2^{n_2} + (\x_{1 k}, \x_{2,d-k})} 
         \right)
- \ell \left( 
       \f{(\x_{1 k}, \x_{2,d-k}) \cap  I_1 I_2^{n_2}}
  { \x_{1 k} I_2^{n_2} + \x_{2,d- k} I_1 I_2^{n_2-1}}
  \right)
  + \sum_{i \geq 2} (-1)^i 
                 h_i(\x_{1 k}, \x_{2,d- k})(1, n_2)\\ \nno
&=& \ell \left(
    \f{R}{ (\x_{1k}, \x_{2,d-k})} 
    \right)
-   \ell \left( 
    \f{I_1^{n_1} I_2^{n_2} + (\x_{1, k}, \x_{2,d-k})} 
     {(\x_{1, k_1}, \x_{2k_2})}
    \right)\\ \nno
&&                 
-   \ell \left( 
    \f{(\x_{1 k_1}, \x_{2 k_2}) \cap  I_1^{n_1} I_2^{n_2}}
      { \x_{1 k}  I_2^{n_2} + \x_{2,d- k} I_1I_2^{n_2-1}}
    \right)
+   \sum_{i \geq 2} (-1)^i 
                 h_i(\x_{1 k}, \x_{2,d- k})(1, n_2)\\ \nno         
&=& e(\x_{1k}, \x_{2,d-k}) 
-    \ell 
     \left( 
     \f{ I_1 I_2^{n}}
             {\x_{1k}   I_2^{n_2}
            + \x_{2,d-k} I_1   I_2^{n_2-1}}
     \right)
+    \sum_{i \geq 2} (-1)^i 
                 h_i(x_{1}, \x_{2,d-k})(1, n_2)\\ \nno
                 &=& \Delta^{k,d-k}P_{\mathcal F}(n_1,n_2)
                 -    \ell 
     \left( 
     \f{ I_1 I_2^{n}}
             {\x_{1k}   I_2^{n_2}
            + \x_{2,d-k} I_1   I_2^{n_2-1}}
     \right)
+    \sum_{i \geq 2} (-1)^i 
                 h_i(x_{1}, \x_{2,d-k})(1, n_2).
                 \eeq \qed
                 
In    \cite[Proposition~2.5]{jay-verma} a version of Huneke's fundamental lemma was obtained for the fiber cone $F_{I_1}(I_2)$ of dimension two where $I_2$ is an $\m$-primary ideal and $I_1 \supseteq I_2$. We obtain a generalization for the fiber cone $F_{I_1}(I_2)$ of  dimension $d \geq 2$ and for any two $\m$-primary ideals $I_1$ and $I_2$ (see Theorem~\ref{fundamental-fibercone}).

 It is well known that for any two $\m$-primary ideals in a local ring $(R, \m)$, the function $H( F_{I_1}(I_2), n) :=  H_{\mathcal F}(1, n)- H_{\mathcal F}(0, n)
 = \ell(I_1^n/I_1 I_2^n)$ is a polynomial say $P( F_{I_1}(I_2), n)$ of degree $d-1$ which can be written in the form 
 \beqn
 P( F_{I_1}(I_2), n)
 = f_{0,I_1} {n + d-1 \choose d-1}
 - f_{1,I_1} {n + d-2 \choose d-2}
 + \cdots
 + (-1)^{d-1}f_{d-1,I_1}.
 \eeqn
 
\blem
  {[Huneke's fundamental Lemma and mixed multiplicities]}
  \label{extension-fiber-jred}
   Let  $I_1$ and $I_2$ be an $\m$-primary ideals in a   be a Cohen-Macaulay ring $(R, \m)$ of dimension $d \geq 2$. 
   \been
   \item
   Let $x_1 \in I_1$ and $\x_{2,d-1} \in I_2$. Assume that $x_1, \x_{2,d-1}$ is   a Rees superficial sequence for $I_1$ and $I_2$. Then
\beq
 \label{extension-fiber-2}
\Delta^{d-1} \left[H( F_{I_1}(I_2), n) \right]    
=   e_{d-1}( I_1 | I_2)  
-     \ell 
       \left(  
       \f{I_1 I_2^{n}}
           {x_{1} I_2^{n} + \x_{2~d-1} I_1 I_2^{n-1}  } 
        \right)
+\sum_{i=2}^{d} (-1)^i h_i(x_{1}, \x_{2~d-1})(1, n).
 \eeq
 
 \item
 Assume that $ \x_{2d} \in I_2$ is   a Rees superficial sequence of $I_1$ and $I_2$. Then
\beq
 \label{extension-fiber-3}
\Delta^{d-1} \left[ H( F_{I_1}(I_2), n)    \right]
=   e( I_2)  
-     \ell 
       \left(  
       \f{I_1 I_2^{n}}
           { x_{2d}I_2^n + \x_{2,d-1} I_1 I_2^{n-1}  } 
        \right)
+\sum_{i=2}^{d} (-1)^i h_i(x_{2d},\x_{2,d-1})(1, n).
 \eeq
 \eeen
\elem
\proof We first prove  (\ref{extension-fiber-2}). Put $n_1 =1$, $n_2 = n$, $k=1$ in 
Corollary~\ref{cor-extension-huneke-d}. The left hand side of 
(\ref{extension-huneke-general}) is 
 \beq
 \label{extension-fiber-d-1}
 && \Delta^{1,d-1} 
  \left[ P_{\mathcal F}(1,n)
  - H_{\mathcal F}(1,n_2) \right]\\ \nno
  &=&    e_{d-1}(I_1|I_2)-     \sum_{i=0}^{d} (-1)^i
     \sum_{j=0}^{i} {1 \choose i-j} {d-1 \choose j}
      H_{\mathcal F}(1 - (i-j),n-j)\\ \nno
      &=& e_{d-1}(I_1|I_2)-    \sum_{i=0}^{d-1} (-1)^i {d-1 \choose i}
        H_{\mathcal F}(1,n-i)
        + \sum_{i=0}^{d-1} (-1)^i {d-1 \choose i}
         H_{\mathcal F}(0,n-i) \\ \nno
         &=& 
          \label{extension-fiber-d-1}e_{d-1}(I_1|I_2)-    \sum_{i=0}^{d-1} (-1)^i{d-1 \choose i} 
H( F_{I_1}(I_2), n-i)\\ \nno
&=& e_{d-1}(I_1|I_2)-   \Delta^{d-1}H( F_{I_1}(I_2), n).
 \eeq
and the right hand side of (\ref{extension-huneke-general}) is 
\beq
 \label{extension-fiber-d-2}
     \ell 
     \left( 
     \f{ I_1 I_2^{n}}
             {x_{1}   I_2^{n}
            + \x_{2,d-1} I_1   I_2^{n-1}}
     \right)  \sum_{i \geq 2} (-1)^i 
                 h_i(x_{1}, \x_{2,d-1})(1, n).
                 \eeq
From (\ref{extension-fiber-d-1}) and 
(\ref{extension-fiber-d-2}) we get (\ref{extension-fiber-2}).

We first prove  (\ref{extension-fiber-3}). Put $n_1 =1$, $n_2 = n$, $k=0$ in 
Corollary~\ref{cor-extension-huneke-d}. The left hand side of 
(\ref{extension-huneke-general}) is 
 \beq
 \label{extension-fiber-d-4}
 && \Delta^{1,d-1} 
  \left[ P_{\mathcal F}(1,n)
  - H_{\mathcal F}(1,n_2) \right]
  = e(I_2)-   \Delta^{d-1}H( F_{I_1}(I_2), n).
 \eeq
and the right hand side of (\ref{extension-huneke-general}) is 
\beq
 \label{extension-fiber-d-5}
     \ell 
     \left( 
     \f{ I_1 I_2^{n}}
             {x_{2d}   I_2^{n}
            + \x_{2,d-1} I_1   I_2^{n-1}}
     \right)  - \sum_{i \geq 2} (-1)^i 
                 h_i(x_{2d}, \x_{2,d-1})(1, n).
                 \eeq
From (\ref{extension-fiber-d-4}) and 
(\ref{extension-fiber-d-5}) we get (\ref{extension-fiber-2}).
\qed
      
     \bnot
We denote
\beqn
        L_{I_1}(I_2) 
&:=& \limn~
  \ell 
  \left( 
   \f{ I_1 I_2^{n}}
   {x_{1} I_2^n + \x_{2,d-1} I_1 I_2^{n-1}  } 
   \right)\\
   L_{I_1,d}(I_2) 
&=& \limn~
  \ell 
  \left( 
   \f{ I_1 I_2^{n}}
   {x_{2d} I_2^n + \x_{2,d-1} I_1 I_2^{n-1}  } 
   \right).
\eeqn
\enot
 
 The following result was proved in \cite[Theorem~2.9]{jay-tony-verma} for the fiber cone $F_{\m}(I)$. We prove it for the fiber cone $F_{I_1}(I_2)$.

\bco
\label{multiplicity-fibercone}
Let $(R, \m)$ be a Cohen-Macaulay local ring of dimension $d \geq 2$. Let $I_1$ and $I_2$ be 
$\m$-primary ideals. 
\been
\item
 Let $x_1 \in I_1$ and $x_{2,d-1} \in I_2$. Assume that $x_1, \x_{2,d-1}$ is   a Rees superficial sequence of $I_1$ and $I_2$. 
Then
\beq
\label{hilb-coef-01}
f_{0, I_1}(I_2) 
 &=& e_{d-1} (I_1 | I_2) -   L_{I_1}(I_2).
\eeq

\item
 Let  $x_{2,d-11}, x_{2,d} \in I_2$ be    a Rees superficial sequence of $I_1$ and $I_2$. Then
 \beq
\label{hilb-coef-02}
f_{0, I_1}(I_2) 
 &=& e_{d-1} (I_2) -   L_{I_1,d}(I_2).
\eeq
\eeen
\eco
   \proof From  Lemma~\ref{vanishing}(\ref{vanishing-two}) we get 
   $ h_i(x_{1}, \x_{2~d-1})(1, n) =   h_i(x_{2d}, \x_{2~d-1})(1, n)= 0$ for all $n \gg 0$. 
   Hence from 
   Lemma~\ref{extension-fiber-jred}
   we have
 \beq
\label{coefficient-fiber} \nno
&&    \limn
\Delta^{d-1}
       H( F_{I_1}(I_2), n-i)    \\
 &=& \label{coefficient-fiber-1}  \limn~ e_{d-1}( I_1 | I_2)  
-    \limn~\ell 
     \left(  
     \f{I_1 I_2^{n}}
       {x_{1} I_2^{n}+ \x_{2,d-1} I_1 I_2^{n-1}  } 
      \right)
= e_{d-1}( I_1 | I_2) 
-   L_{I_1}(I_2)\\
 &=& \label{coefficient-fiber-2}  \limn~ e( I_2)  
-    \limn~\ell 
     \left(  
     \f{I_1 I_2^{n}}
       {x_{2d} I_2^{n}+ \x_{2,d-1} I_1 I_2^{n-1}  } 
      \right)
= e(I_2) 
-   L_{I_1,d}(I_2).
\eeq
\qed

We now state and prove the Huneke's fundamental lemma for fiber cones.

      \bt 
      {[Huneke's fundamental lemma for the fiber cone]}
      \label{fundamental-fibercone}
            With the assumptions as in Lemma~\ref{extension-fiber-jred}(1) we have
      \beqn
      \Delta^{d-1}
\left[        P( F_{I_1}(I_2), n) -  H( F_{I_1}(I_2), n) \right]
=  \ell 
       \left(  
       \f{I_1 I_2^{n}}
           {x_{1} I_2^{n} + \x_{2~d-1} I_1 I_2^{n-1}  } 
        \right) - L_{I_1}(I_2)
- \sum_{j=2}^{d} (-1)^j h_j(x_{1}, \x_{2~d-1})(1, n).
 \eeqn
 \et
    \proof The proof follows from Lemma~\ref{extension-fiber-jred} and Corollary~\ref{multiplicity-fibercone}.

 \section{Hilbert coefficients of the fiber cone}

     We now completely relate the   Hilbert coefficients $f_{i, I_1}(I_2)$ with the homology modules $h_i(x_{1}, \x_{2~d-1})(1, n)$.
  We are able to extend the results on the coefficients of the associated graded ring in \cite{huc} to the fibercone.

Recall that 
$H( F_{I_1}(I_2), n) 
$ is a polynomial in $n$ which we will denote by 
$P( F_{I_1}(I_2), n) 
$ and this polynomial can be written in the form
\beqn
P( F_{I_1}(I_2), n)
= \sum_{i=0}^{d-1} (-1)^i 
f_{i, I_1}
{n + d-1-i \choose d-1-i}.
\eeqn


The following result is well known. We state it for the state of completeness.

\blem
\label{hilb-poly-zero-1}
Let $I_1$  and $I_2$ be be an $\m$-primary ideals in a Cohen-Macaulay ring $(R, \m)$. Assume that there exists a sequence $\x_{2d} \in I_2$ such that 
 $\x_{d-1}^{o}$ (resp. $\x^{*}_{d-1}$) is a superficial sequence in $F_{I_1}(I_2)$ (resp. $G(I_2)$). Then
for all $0 \leq i \leq d-1$,
\been
\item
$
 \displaystyle{      \Delta^{d-1-i} P( F_{I_1}(I_2), -1)
=   (-1)^i f_{i, I_1}(I_2) .}
$

\item
$
 \displaystyle{      \Delta^{d-1-i} P( F_{I_1}(I_2), 0)
=    \sum_{j=0}^i(-1)^j f_{j, I_1}(I_2) .}
$
\eeen
\elem
\proof The proof can be easily verfied for $d=1$ and for $d>1$ and $i=d-1$. 

Now let $d > 1$ and let $i < d-1$.   Let $-$ denote the image in $R/x_{21}$. By our assumption of $x_{21}$, 
\beqn
 \dim~R/x_1 &=& \dim~R -1\\
 \dim~F_{  \olin{I_1}}(    {\olin{I_2}}) &=& \dim~F_{I_1}(I_2) -1\\
 \Delta [P( F_{I_1}(I_2), n)]
&=& P( F_{   {\olin{I_1}}}    (  {\olin{I_2}}), n)
\\
f_{i, I_1}(I_2) &=& f_{i,  \olin{I_1}}( \olin{I_2}), \hspace{.5in} \mbox{ for all } i = 0, \ldots, d-2.  
\eeqn
Therefore for all  $i=0, \ldots, d-2$, 
\beqn
\Delta^{d-1-i} P_{I_1}(I_2,-1) 
=   \Delta^{d-2-i} P_{{\olin I_1}}({\olin I_2},-1) 
=  (-1)^i  f_{i,  \olin{I_1}}( \olin{I_2})
=  (-1)^i f_{i, I_1}(I_2) .
\eeqn
This proves(1). The proof of (2) is similar.
\qed

\bco
\cite[Lemma~2.7]{huc}
\label{cor-lemma-2.7-1}
Let $f: \Z \rar \Z$ be a function such that $f(n) = 0$ for $n \gg 0$. Then
\been
\item
$\displaystyle
{\sum_{n=i-1}^{\infty} {n \choose i-1} \Delta^j f(n)
= (-1)^i \Delta^{j-i} f(-1)
}.
$

\item
 $\displaystyle
{
\sum_{n=i}^{\infty} {n-1 \choose i-1} \Delta^j f(n)
= (-1)^i \Delta^{j-i} f(0)}.
$
\eeen
\eco
\proof The proof is similar to the proof of Lemma~2.7 of \cite{huc}.
\qed



We now extend Proposition~2.9 of \cite{huc} to the Hilbert coefficients of the fiber cone.

 \bp 
 \label{fiber-c}
 Let $(R, \m)$ be a Cohen-Macaulay local ring. Let $I_1$ and $I_2$ be $\m$-primary ideals of $R$. Then
 \been
 \item
\label{f-i-1}
$\displaystyle
{f_{i, I_1}(I_2) 
= \sum_{n \geq i-1}{n \choose i-1} \Delta^{d-1} 
\left[  P( F_{I_1}(I_2), n) - H( F_{I_1}(I_2), n)  \right]}
.
$

\item
$\displaystyle
{ (-1)^{i}\sum_{j=0}^i (-1)^j f_{j, I_1}(I_2) 
 + (-1)^{i+1} \ell \left( \f{R}{I_1} \right)
 =
\sum_{n \geq i}{n-1 \choose i-1} \Delta^{d-1} 
\left[  P( F_{I_1}(I_2), n) - H( F_{I_1}(I_2), n)  \right]
}.$
\eeen
\ep
\proof Put $j = d-1$ and $f(n) =P( F_{I_1}(I_2), n) - H( F_{I_1}(I_2), n)$ in Corollary~\ref{cor-lemma-2.7-1}. 
Then
\beqn
            \sum_{n=i-1}^{\infty} {n \choose i-1} 
            \Delta^{d-1} 
            \left[  P( F_{I_1}(I_2), n) - H( F_{I_1}(I_2), n) \right]
&=&  (-1)^i \Delta^{d-1-i} 
            \left[ P( F_{I_1}(I_2), -1) - H( F_{I_1}(I_2), -1) \right] \\
&=&  (-1)^i \Delta^{d-1-i} 
            \left[ P( F_{I_1}(I_2), -1)\right]  \\
&=&   f_{i, I_1}(I_2) \hphantom{oneinch} \mbox{[by Lemma~\ref{hilb-poly-zero-1}(1)]}.
\eeqn
Similarly,
\beqn
           \sum_{n=i}^{\infty} {n-1 \choose i-1} 
           \Delta^{d-1} 
           \left[  P( F_{I_1}(I_2), n) - H( F_{I_1}(I_2), n) \right]
&=& (-1)^i \Delta^{d-1-i} 
          \left[ P( F_{I_1}(I_2), 0)  - H( F_{I_1}(I_2), 0) \right] \\
&=&  (-1)^{i+1}\sum_{j=0}^i (-1)^j f_{j, I_1}(I_2)
- \ell \left( \f{R}{I_1} \right).
\eeqn
\qed


We now describe the Hilbert coefficeints of the fiber cones in terms of the homology modules \\
$h_i(x_{1}, \x_{2~d-1})(1, *)$. 
\bt
\label{thm-f-i}
Let $I_1$ and $I_2$ be $\m$-primary ideals in a Cohen-Macaulay local ring $(R,\m)$. Then 
\label{f-i-2}
\beqn
f_{1, I_1}(I_2) 
&=& f_{0, I_1}(I_2)  - \ell \left( \f{R}{I_1} \right) \\
&&+  \sum_{n \geq 1} \left[ 
      \ell 
      \left( \f{I_1 I_2^n}
                {x_1 I_2^n +  \x_{2,d-1} I_1 I_2^{n-1})} 
                \right)
-        L_{I_1}(I_2) \right]
 -    \sum_{j \geq 2} (-1)^{j} 
                        h_j(x_{1}, \x_{2~d-1})(1, *)\\
        f_{i, I_1}(I_2) 
&=&    \sum_{n \geq i-1}
        {n \choose i-1}
         \left[ \ell 
      \left( \f{I_1 I_2^n}
                {x_1 I_2^n +  \x_{2,d-1} I_1 I_2^{n-1})} 
                \right)
-        L_{I_1}(I_2)  -    \sum_{j \geq 2} (-1)^{j} 
                        h_j(x_{1}, \x_{2~d-1})(1, n)\right], i \geq 2.
 \eeqn
\et
\proof By Proposition~\ref{fiber-c} we get
\beqn
    f_{1, I_1}(I_2) 
&=& \sum_{n \geq 0} \Delta^{d-1} 
     \left[  P( F_{I_1}(I_2), n) - H( F_{I_1}(I_2), n)  \right]\\
&=& \Delta^{d-1} 
    \left[  P( F_{I_1}(I_2), 0) - H( F_{I_1}(I_2), 0)  \right]
+   \sum_{n \geq 1} \Delta^{d-1} 
    \left[  P( F_{I_1}(I_2), n) - H( F_{I_1}(I_2), n)  \right]\\
&=&  f_{0, I_1}(I_2)  - \ell \left( \f{R}{I_1} \right) \\
&&+  \sum_{n \geq 1} \left[ 
      \ell 
      \left( \f{I_1 I_2^n}
                {x_1 I_2^n +  \x_{2,d-1} I_1 I_2^{n-1})} 
                \right)
-        L_{I_1}(I_2) \right]
 -    \sum_{n \geq 1} \sum_{i \geq 2} (-1)^{i} 
                        h_i(x_{1}, \x_{2~d-1})(1, n)\\
                        &=&  f_{0, I_1}(I_2)  - \ell \left( \f{R}{I_1} 
                        \right) \\
                       &&+  \sum_{n \geq 1} \left[ 
      \ell 
      \left( \f{I_1 I_2^n}
                {x_1 I_2^n +  \x_{2,d-1} I_1 I_2^{n-1})} 
                \right)
-        L_{I_1}(I_2) \right]
 -    \sum_{i \geq 2} (-1)^{i} 
                        h_i(x_{1}, \x_{2~d-1})(1, *).
 \eeqn
The third equality follows from Lemma~\ref{fundamental-fibercone}.

For $j \geq 1$ it follows from  Lemma~\ref{fundamental-fibercone} and Proposition~\ref{fiber-c}.
\qed

We now describe the relation between the different Hilbert coefficeints of the fiber cones and the homology modules $h_i(x_{1}, \x_{2~d-1})(1, *)$. 
\bt
Let $I_1$ and $I_2$ be $\m$-primary ideals in a Cohen-Macaulay local ring $(R,\m)$.  Let $i \geq 2$. Then 
\label{f-i-2}
\beqn
 && (-1)^i\sum_{j=0}^{i}
(-1)^{j}  f_{j, I_1}(I_2)  
+ (-1)^{i+1}\ell \left( \f{R}{I_1} \right)\\
&=&  \sum_{n \geq i-1}
        {n-1 \choose i-1}
         \left[ \ell 
      \left( \f{I_1 I_2^n}
                {x_1 I_2^n +  \x_{2,d-1} I_1 I_2^{n-1})} 
                \right)
-        L_{I_1}(I_2)  -    \sum_{k \geq 2} (-1)^{k} 
                        h_l(x_{1}, \x_{2~d-1})(1, n)\right].
                        \eeqn
                        \et
                        \proof  One can show that 
                        \beqn
\sum_{n \geq i}{n-1 \choose i-1} \Delta^{d-1} 
\left[  P( F_{I_1}(I_2), n) - H( F_{I_1}(I_2), n)  \right]
= (-1)^{i}\sum_{j=0}^i (-1)^j f_{j, I_1}(I_2) .
\eeqn
Now use  Lemma~\ref{fundamental-fibercone}.
                        \qed

  \bco
  \label{cor-bound}
  Let $I_1$ and $I_2$ be $\m$-primary ideals in a Cohen-Macaulay local ring $(R, \m)$.
\beq
\label{cor-bound-f-1}
        f_{1, I_1}(I_2) 
&\leq&  f_{0, I_1}(I_2)  - \ell \left( \f{R}{I_1} \right)
+ \sum_{n \geq 1} \left[ 
      \ell 
      \left( \f{I_1 I_2^n}
                {x_1 I_2^n +  \x_{2,d-1} I_1 I_2^{n-1})} 
                \right)
-        L_{I_1}(I_2) \right]
.
\eeq
If equality holds in (\ref{cor-bound-f-1}), then
for all $j \geq 2$, 
\beq
\label{bound-f-2}
 f_{j, I_1}(I_2) 
&=&    \sum_{n \geq j-1}
        {n \choose j-1}
         \left[ \ell 
      \left( \f{I_1 I_2^n}
                {x_1 I_2^n +  \x_{2,d-1} I_1 I_2^{n-1})} 
                \right)
-        L_{I_1}(I_2)  \right] \\
\label{bound-f-3}
\sum_{j=0}^{i}
(-1)^{j}  f_{j, I_1}(I_2)  
- \ell \left( \f{R}{I_1} \right)
&=&  \sum_{n \geq i-1}
        {n-1 \choose i-1}
         \left[ \ell 
      \left( \f{I_1 I_2^n}
                {x_1 I_2^n +  \x_{2,d-1} I_1 I_2^{n-1})} 
                \right)
-        L_{I_1}(I_2) \right]
\eeq

\eco
\proof  By Lemma~\ref{rigidity}, 
\beqn
 \sum_{i \geq 2} (-1)^{i} h_i(x_{1}, \x_{2 k_2})(1, *) \geq 0.
\eeqn
Hence by Theorem~\ref{f-i-1}, we get (\ref{bound-f-1}). If equality holds in (\ref{bound-f-1}), then 
\beqn
 \sum_{i \geq 2} (-1)^{k-2} h_i(x_{1}, \x_{2 k_2})(1, *) = 0.
\eeqn
then by Lemma~\ref{rigidity},  $h_i(x_{1}, \x_{2 k_2})(1, n) =0$ for all $j \geq 2$ and for all $n$. Substituting $h_i(x_{1}, \x_{2 k_2})(1, n) =0$ for all 
$i \geq 2$ and for all $n$ in Theorem~\ref{f-i-1}
and  Theorem~\ref{thm-f-i} we get  (\ref{bound-f-2}) and (\ref{bound-f-3}).
\qed

\section{Vanishing and depth of associated graded ring and fiber cone}

In this section we are interested in interlinking the depth of fiber cone, the depth of the associated graded ring. This has been studied for ideals of minima mixed multiplicity \cite{clare-verma1}, ideals of almost minimal mixed multiplicity \cite{clare-verma2} and for ideals of almost minimal multiplicity \cite{jay-verma}.

The vanishing of the complex
$C({x}_{1}, {\x}_{2k}, {\mathcal F}, (1, n_2))$ plays an important role. Earlier results in \cite{clare-verma1}, \cite{clare-verma2} and \cite{jay-verma}) can easily be recovered from this (see \cite{clare1} and \cite{clare2}).

\bnot
{\em 
Let $I_1,I_2$ be ideals in a ring $(R, \m)$. For an element $x \in I_2$, let 
$x^{*}$ (resp. $x^{o}$) denote the residue class in $I_2/ I_2^2$ (resp. $I_2/ 
I_1 I_2$).
}
\enot

We have the exact sequence of complexes:
\beq
\label{koszul-fiber}
     0
\rar K({\x_{2k}^{o}}, F_{I_1}(I_2))(n)
\rar C( {\x}_{2k}, {\mathcal F}, (1, n))
\rar  C( {\x}_{2k}, {\mathcal F}, (0, n))
\rar 0
\eeq
where  $K({\x_{2k}}^o, F_{I_1}(I_2))(n)$ is the Koszul complex of the fiber cone with respect to the sequence $\x_{2k}^o = x_{21}^o, \ldots, x_{2k}^o$.
Hence we have the   sequence:
\beq
\label{exact-depth-2}
\begin{array}{lrrl}
     & \cdots
\rar & H_{i+1}(C( \x_{2k_2}, {\mathcal F}, (1, n)))
\rar & H_{i+1}(C( {x}_{1}, {\x}_{2k}, {\mathcal F}, (1, n)))\\
\rar & H_{i}  (K( {\x_{2k_2}}, F_{I_1}(I_2)))(n)
\rar & H_{i}  (C( {\x}_{2k_2}, {\mathcal F}, (1, n)))
\rar & H_{i}  (C( {x}_{1}, {\x}_{2k}, {\mathcal F}, (1, n))
\end{array}.
\eeq

Recall that we also have the  exact sequence:
\beqn
     0
\rar {{C(           {\x}_{2k},  {\mathcal F}, (1, n))}}
\rar {{C( {x}_{1},  {\x}_{2k}, {\mathcal F}, (1, n))}}
\rar {{C(           {\x}_{2k},  {\mathcal F}, (0, n))}}
\rar 0.
\eeqn
and  the   sequence:
\beq
\label{exact-depth-0}
          \cdots
\rar &  H_{i+1}( {{C({x}_{1}, {\x}_{2k}, {\mathcal F}, (1, n))}} )
\rar &  H_{i}  ({{C(          {\x}_{2k},  {\mathcal F}, (0, n))}})\\ \nno
\rar    H_{i}  ({{C(          {\x}_{2k},   {\mathcal F}, (1, n))}})
\rar&   H_{i}  ({{C({x}_{1},  {\x}_{2k},   {\mathcal F}, (1, n))}} )
\rar&   \cdots
\eeq

\bnot
{\em 
 Let $\gamma(\x^{*},G(I_2)) = \grade~(x^{*},G(I_2))$
  (resp.   $\gamma( \x^{o}F_{I_1}(I_2)) = \grade~(\x^{o}, F_{I_1}(I_2)$
 denote the depths of the ideals 
 $(\x_{2k}^{*})$
  (resp.  $(\x_{2k}^o)$
in $G(I_2)$
(resp.
$F_{I_1}(I_2)$).
}
\enot

\bt
\label{depth-fiber}
Let $I_1$ and $I_2$ be $\m$-primary ideals in a Cohen-Macaulay local ring $(R,\m)$  of dimension $d \geq 2$. Let $x_{1} \in I_1$ and $ \x_{2,d-1} \in I_2$ be such that $x_1, \x_{2,d-1}$ is  a superficial sequence for  $1$ copy of $I_1$ and $k$ copies of $I_2$. 
Let
\beqn
i = \min \{j | H_{d-j}( {{C({x}_{1}, {\x}_{2k}, {\mathcal F}, (1, n))}} )\not =0  \mbox{ for some } n\}.
\eeqn
Let $g =  \gamma(\x_{2k}^{*}, G(I_2))$.
\been
\item
\label{koszul-fiber-i}
If $g \leq i-1$, then  $\gamma(x_{2k}^{o}F_{I_1}(I_2)) \geq g$.

\item
\label{koszul-fiber-ii}
If $g \geq i$, then  $ \gamma(\x_{k}^{*}, F_{I_1}(I_2)) \geq i$.
 \eeen
\et
\proof 
Suppose $g \leq i-1$, by \cite[Proposition~3.3]{tom-huc},
\beq
\label{depth-3}
H_{d-1-g+a} C(          {\x}_{2k},   {\mathcal F}, (o, n))
= 0,  \mbox{for all } a \geq 1 \hphantom{oneinch} \mbox{ [\cite{anna-thesis}, \cite{tom-huc}]}.
\eeq 
By assumption
\beqn
H_{k-i+a} C(x_1,    {\x}_{2k},   {\mathcal F}, (1, n)) = 0, 
\hspace{.5in} \mbox{ for all } 
a \geq 2.
\eeqn
  Substituting in (\ref{exact-depth-0}) we get
 \beq
 \label{depth-4}
H_{k-g+a} C(          {\x}_{2k},   {\mathcal F}, (0, n)) = 0, 
\hspace{.5in} \mbox{for all } 
a \geq 1.  
\eeq
Substituting (\ref{depth-3}) and (\ref{depth-4}) in (\ref{exact-depth-2}) we get
\beqn
H_{k-g+a}  (K( {\x_{2k}^{o}}, F_{I_1}(I_2)))(n)
\hspace{.5in} \mbox{for all } 
a \geq 1. 
\eeqn
Hence $\gamma(x_{2k}^{o}, F_{I_1}(I_2)) \geq g$.

Similarly if $g >i$, then one can show that 
 \beqn
H_{k-i+a}  (K( {\x_{2k}}, F_{I_1}(I_2)))(n) =0
\hspace{.5in} \mbox{for all } 
a \geq 1. 
\eeqn
Hence $\gamma(x_{2k}^{o}, F_{I_1}(I_2)) \geq i$.
\qed

\bco
Let $I_1$ and $I_2$ be $\m$-primary ideals in a Cohen-Macaulay local ring $(R, \m)$. Let $x_1 \in I_1$ and $\x_{d-1} \in I_2$. Let $x_1, \x_{2,d-1}$ be a superficial sequence for $I_1$ and $I_2$. Suppose
\beqn
d-2 = \min \{j | H_{d-j}( {{C({x}_{1}, {\x}_{2k}, {\mathcal F}, (1, n))}} )\not =0  \mbox{ for some } n\}.
\eeqn
If  $\gamma(\x_{2k}^{*}, G(I_2)) \geq d-1$, then $\gamma(\x_{2k}^{o}F_{I_1}(I_2)) \geq d-1$.
\eco
\proof Put $k=d$ in Theorem~\ref{depth-fiber}.

\bt 
Let $I_1$ and $I_2$ be $\m$-primary ideals in a Cohen-Macaulay local ring $(R, \m)$. If
\beq
\label{bound-f-1}
        f_{1, I_1}(I_2) 
&=&  f_{0, I_1}(I_2)  - \ell \left( \f{R}{I_1} \right)
+  \sum_{n \geq 1} \left[ 
      \ell 
      \left( \f{I_1 I_2^n}
                {x_1 I_2^n +  \x_{2,d-1} I_1 I_2^{n-1})} 
                \right)
-        L_{I_1}(I_2) \right]
\eeq
and $\gamma( \x_{2,d-1}^{*}, G(I_2)) \geq d-1$
then $\gamma(\x_{2,d-1}^{o} F_{I_1}(I_2)) \geq d-1$.
\et
\proof By Corollary~\ref{cor-bound},
\beqn
 \min \{j | H_{d-j}( {{C({x}_{1}, {\x}_{2k}, {\mathcal F}, (1, n))}} )\not =0  \mbox{ for some } n\} \geq d-1.
\eeqn
Now apply Theorem~\ref{depth-fiber}(2).
\qed

   \end{document}